\documentclass[12pt]{article}

\usepackage[left=1.5in, right=1in, top=1in, bottom=1in, dvips,pdftex]{geometry}
\usepackage{color, setspace, graphicx, verbatim, mathrsfs, amsmath, amsthm, amssymb, amsfonts,  multirow, moreverb}
\newtheorem{defn}{Definition}[section]

\newtheorem{thm}[defn]{Theorem}

\theoremstyle{remark}

\numberwithin{equation}{section} \numberwithin{figure}{subsection}

\DeclareMathOperator{\airy}{Ai}

\def\ra{\rightarrow}
\def\iy{\infty}

\def\be{\begin{equation}}
\def\ee{\end{equation}}

\newcommand{\bP}{\mathbb{P}}

\begin{document}

\title{\textbf{Finite $n$  Largest Eigenvalue Probability Distribution Function
 of Gaussian Ensembles }}

\author{Leonard N.~Choup \\ Department of Mathematical Sciences \\ University Alabama in Huntsville\\
 Huntsville, AL 35899, USA \\
 email:  \texttt{Leonard.Choup@uah.edu}
} \maketitle

\begin{center} \textbf{Abstract} \end{center}

\begin{small}
In this paper we focus on the finite $n$ probability distribution
function of the largest eigenvalue in the classical Gaussian
Ensemble of $n\times n$ matrices (GE$_n$). We derive the finite $n$ largest eigenvalue probability distribution function for the Gaussian Orthogonal and Symplectic Ensembles  and also
prove an Edgeworth type Theorem for the largest eigenvalue probability distribution
function of Gaussian Symplectic Ensemble. The correction
terms to the limiting probability distribution are expressed in
terms of the same Painlev$\acute{e}$ II functions appearing in the
Tracy-Widom distribution.

\end{small}

\section{Introduction}\label{introduction}
In applications of the limiting probability distributions laws from
Random Matrix Theory (e.g. \cite{Deif1}, \cite{John1}) it is
important to have an estimate on the convergence rates, if possible
have a control on this rate of convergence. For recent reviews of
applications of these distributions we refer the reader to
\cite{Deif2,Deif3, Dien1, Johansson, Trac8}. In our desire to
control the rate of convergence of the probability distribution of
the largest eigenvalue from the Gaussian Orthogonal Ensemble $GOE_{n}$, Gaussian Unitary Ensemble $GUE_{n}$ and
Gaussian Symplectic Ensembles $GSE_{n}$, we introduce a fine tuning constant $c$ in the
scaling of the desired eigenvalue. The finite large $n$ expansion is
therefore a function of $c$. We use this constant to fine-tune the
convergence rate. This work was done for the $GUE_{n}$ and $GOE_{n}$. In completing the same work for the $GSE_{n}$, we decided 
to find a closed formula (opposed to a large $n$ asymptotic formula) whenever possible for each function appearing in the the final expression for these probabilities distribution functions. This approach has the advantages that the final results is a finite $n$ representation of the distributions functions, we only need to perform the large $n$ expansion once for each function to have an Edgeworth type expansion for the desired probability.

The Gaussian
$\beta$–-ensembles are probability spaces on $n$-tuples of random
variables $\{\lambda_{1}, \ldots , \lambda_{n}\}$ (think of them as
eigenvalues of a randomly chosen matrix from the ensemble.) with
probability density that the variables lie in an infinitesimal
intervals about the points $x_{1}, \ldots , x_{n}$ is
\begin{equation}\label{jpdfeig}
\bP_{n\beta}(x_{1},\cdots,x_{n})\; = \textrm{C}_{n\beta}\,
\textrm{exp}\left(-\frac{\beta}{2}\, \sum_{1}^{n}x_{j}^{2}\right)\,
\prod_{j<k}|x_{j}-x_{k}|^{\beta},
\end{equation}
with
\begin{equation}
-\infty < x_{i} < \infty, \quad \textrm{for} \; i=1, \cdots, n.
\end{equation}
Here $\textrm{C}_{n\beta}$ is the normalizing constant such that the
total integral over the $x_{i}$'s is one. When
$\beta=1$ we have the GOE$_{n}$, when $\beta=2$ we have the GUE$_{n}$ and
when $\beta=4$ we have the GSE$_{n}$. We denote the largest eigenvalue by
$\lambda_{Max}^{\beta}$, and by
\begin{equation}\label{p.d.f}
F _{n,\beta}(t)= \bP(\lambda_{\textrm{max}}^{\beta} \leq t)
\end{equation}
his probability distribution function.\\
When $\beta=2$, the harmonic oscillator wave functions (see
\cite{Hoch1}, or \cite{Szeg1} for a complete definition)
\begin{equation*}
\varphi_{k}(x)= {1\over (2^{k}k! \sqrt{\pi})^{1/2}} \, H_k(x)\,
e^{-x^2/2}\quad k=0,\,1,\,2,\ldots
\end{equation*}
  play an important role.  We also
have the Hermite kernel
\begin{equation}\label{hermite kernel}
K_{n,2}(x,y)=\sum_{k=0}^{n-1}\varphi_{k}(x)\varphi_{k}(y) =
\sqrt{\frac{n}{2}} \quad \frac{\varphi_{n}(x)\varphi_{n-1}(y)
 -  \varphi_{n}(y)\varphi_{n-1}(x)}{x-y},
\end{equation}
which is the kernel of the integral operator $K_{n,2}$ acting on
$L^{2}(t,\iy)$ , with resolvent kernel
\begin{equation}\label{resolvent kernel}
R_{n,2}(x,y;t)\>=\> (I-K_{n,2})^{-1}\cdot K_{n,2}(x,y).
\end{equation}
Note here that the dot denotes operator multiplication.
 We have the following representation of
\eqref{jpdfeig}, (see for example, \cite{Meht1} for a derivation of this result)
\begin{equation*}
\bP_{n2}(x_{1},\cdots,x_{n})\>=\> \det(K_{n,2}(x_{i},x_{j}))_{1\leq
i\,j\,\leq n}.
\end{equation*}
Following Tracy and Widom in  \cite{Trac3, Trac7, Trac2, Trac1,
Trac4, Trac8, Trac5}, we
define
\begin{equation}\label{varphi}
\varphi(x)\>=\>
\biggl(\frac{n}{2}\biggr)^{\frac{1}{4}}\varphi_{n}(x),\quad \quad
\psi(x)\>=\>
\biggl(\frac{n}{2}\biggr)^{\frac{1}{4}}\varphi_{n-1}(x),
\end{equation}
by $\varepsilon$ the integral operator with kernel
\begin{equation}\label{varepsilon}
\varepsilon_{t}(x)\>\>=\>\> \frac{1}{2}\mathrm{sgn}(x-t),
\end{equation}
$D$ the differentiation with respect to the independent variable,
\begin{equation}\label{Qn}
Q_{n,i}(x;t)\, =\, (\,(I-K_{n,2})^{-1}\, ,\, x^{i}\varphi)
\end{equation}
and
\begin{equation}\label{Pn}
P_{n,i}(x;t)\, =\, (\,(I-K_{n})^{-1}\, ,\, x^{i}\psi).
\end{equation}
We introduce the following quantities
\begin{equation}\label{qn}
q_{n,i}(t)\>=\> Q_{n,i}(t;t), \>\>\> p_{n,i}(t)\>=\> P_{n,i}(t;t)
\end{equation}
\begin{equation}\label{un}
u_{n,i}(t)\,=\, (Q_{n,i},\varphi),\quad v_{n,i}(t)\,=\,
(P_{n,i},\varphi),\quad
\end{equation}
\begin{equation}
\tilde{v}_{n,i}(t)=(Q_{n,i},\psi),\quad \mathrm{and} \quad
w_{n,i}(t)\,= \, (P_{n,i},\psi).
\end{equation}
 Here $(\,\cdot \,,\cdot \,)$ denotes
the inner product on $L^{2}(t,\infty)$. In our notation, the
subscript without the $n$ represents the scaled limit of that
quantity when $n$ goes to infinity, and  we dropped the second
subscript $i$ when it is zero.\\  If $\airy$ is the Airy function,
the kernel $K_{n,2}(x,y)$ then scales\footnote{as $n\rightarrow \iy$ in the change $x=\sqrt{2(n+c)} +2^{-\frac{1}{2}}n^{-\frac{1}{6}}X$ and 
$y=\sqrt{(n+c)} +2^{-\frac{1}{2}}n^{-\frac{1}{6}}Y$} to the  Airy kernel
\begin{equation}\label{Airy kernel}
K_{\airy}(X,Y)\>=\>
\frac{\airy(X)\,\airy'(Y)\>-\>\airy(Y)\,\airy'(X)}{X-Y}.
\end{equation}
Our conventions are as follows:
\begin{equation}\label{Q}
Q_{i}(x;s)\, =\, (\,(I-K_{\airy})^{-1}\, ,\, x^{i}\airy),\quad
Q_{0}(x;s)\,=\,Q(x;s),
\end{equation}
\begin{equation}\label{P}
P_{i}(x;s)\, =\, (\,(I-K_{\airy})^{-1}\, ,\, x^{i}\airy'),\quad
P_{0}(x;s)\,=\,P(x;s),
\end{equation}
\begin{equation}\label{q}
q_{i}(s)\>=\> Q_{i}(s;s),\quad q_{0}(s)=q(s) ,\>\>\> p_{i}(s)\>=\>
P_{i}(s;s),\quad p_{0}(s)=p(s),
\end{equation}
\begin{equation}\label{u}
u_{i}(s)\,=\, (Q_{i},\airy),\quad u_{0}(s)=u(s),\quad v_{i}(s)\,=\,
(P_{i},\airy),\quad v_{0}(s)=v(s),
\end{equation}
\begin{equation}
\tilde{v}_{i}(s)=(Q_{i},\airy'),\quad
\tilde{v}_{0}(s)=\tilde{v}(s),\quad  w_{i}(t)\,= \,
(P_{i},\airy'),\quad \textrm{and}\quad w_{0}(s)=w(s).
\end{equation}
 Here $(\,\cdot \,,\cdot \,)$  denotes
the inner product on $L^{2}(s,\infty)$  and $i=0,1,2,\cdots$ \\
We also note that $q(s)$ is the Haskins-Macleod solution to the
Pailev$\acute{e}$ II equation $q''(s)=sq(s)+2q^{3}(s)$ with the
boundary condition
$q(s)\sim \mathrm{Ai}(s)$ as $s \rightarrow \infty$. \\
We use the subscript $n$ for unscaled quantities only. 
\begin{equation}
\mathcal{R}_{n,1}:=\int_{-\infty}^{t}R_{n}(x,t;t)\mathrm{d}x,\quad
\mathcal{P}_{n,1}:= \int_{-\infty}^{t}P_{n}(x;t)\mathrm{d}x,\quad
\mathcal{Q}_{n,1}:= \int_{-\infty}^{t}Q_{n}(x;t)\mathrm{d}x,
\end{equation}
and
\begin{equation*}
\mathcal{R}_{n,4}(t):=\int_{-\infty}^{\infty}\varepsilon_{t}(x)R_{n}(x,t;t)\mathrm{d}x,\quad
\mathcal{P}_{n,4}(t):=
\int_{-\infty}^{\infty}\varepsilon_{t}(x)P_{n}(x;t)\mathrm{d}x,
\end{equation*}
\begin{equation}
\mathcal{Q}_{n,4}(t):=
\int_{-\infty}^{\infty}\varepsilon_{t}(x)Q_{n}(x;t)\mathrm{d}x.
\end{equation}
The epsilon quantities are
\begin{equation}\label{Q epsilon}
Q_{n,\varepsilon}(x;t)\> =\> \bigl(\,(I-K_{n})^{-1}(x,y)\, ,\,
\varepsilon \varphi(y)\bigr), \quad
q_{n,\varepsilon}(t)\>=\>Q_{n,\varepsilon}(t;t)
\end{equation}

\begin{equation}\label{u epsilon}
u_{n,\varepsilon}(t)\,=\,
\bigl(Q_{n,\varepsilon}(x;t),\varphi(x)\bigr),\quad
\tilde{v}_{n,\varepsilon}(t)\,=\, \bigl(Q_{n,\varepsilon}(x;t),
\psi(x)\bigr),
\end{equation}
where $(\,\cdot \,,\cdot \,)$  denotes the inner product on
$L^{2}(t,\infty)$. \\

The GOE$_{n}$ and GSE$_{n}$ analogue of \eqref{GUE Edgeworth} in
Theorem \ref{Fn2} bellow will follow from representations (equations
(40) and (41) of \cite{Trac2}.)
\begin{equation}\label{f_{n,1}}
F_{n,1}(t)^{2}
       = F_{n,2}(t)\cdot
       \left( \bigl(1-\tilde{v}
       _{n,\varepsilon}(t)\bigr)\bigl(1-\frac{1}{2}\mathcal{R}_{n,1}(t)\bigr)-\frac{1}{2}\bigl(q_{n,\varepsilon}(t)
       -c_{\varphi}\bigr)\mathcal{P}_{n,1}(t)\right)
\end{equation}
and

\begin{equation}\label{f_{n,4}}
F_{n,4}(t/\sqrt{2})^{2}= F_{n,2}(t)\cdot\left(
       \bigl(1-\tilde{v}_{n,\varepsilon}(t)\bigr)\bigl(1+\frac{1}{2}\mathcal{R}_{n,4}(t)\bigr)+
       \frac{1}{2}q_{n,\varepsilon}(t)\,\mathcal{P}_{n,4}(t)\right).
\end{equation}

To complete the work, we will use the close form  of
$\mathcal{R}_{n,1},\>\mathcal{P}_{n,1},\>\mathcal{R}_{n,4},\>\mathcal{P}_{n,4},\>\tilde{v}_{n,\varepsilon},\>\mathrm{and}\>q_{n,\varepsilon}$ 
derived in \cite{Choup4} to give a finite $n$ representation of $F_{n,1}$ and $F_{n,4}$, then the corresponding large $n$ expansion to find an Edgeworth type expansion for $F_{n,4}$ as outlined in \cite{Choup3}.
We will also need the following results, 
\begin{thm}\cite{Choup1}\label{Fn2}
If we set
\begin{equation}\label{scale for GUE}
\tau( s ) =  (2(n+c)
)^{\frac{1}{2}}+2^{-\frac{1}{2}}n^{-\frac{1}{6}}\,s\>\>\>
\textrm{and}
\end{equation}
\begin{equation}\label{second term}
E_{c,2}:=E_{c,2}(s)=2w_{1}-3u_{2}+ (-20c^2 +3)v_{0}  +
u_{1}v_{0}-u_{0}v_{1} +u_{0}v_{0}^{2}-u_{0}^{2}w_{0}.
\end{equation}
Then as $n\ra\iy$
\begin{equation}\label{GUE Edgeworth}
F_{n,2}(\tau(s))= F_{2}(s)\left\{ 1 + c \, u_{0}(s) \, n^{-\frac{1}{3}}
-\frac{1}{20}E_{c,2}(s)\, n^{-\frac{2}{3}}\right\} + O(n^{-1})
\end{equation}
uniformly in $s$, and
\begin{equation}\label{TW}
F_{2}(s)\>=\> \lim_{n\ra\iy}F_{n,2}(t)\>=\>
\exp\left(-\int_{s}^{\iy}(x-s)q(x)^2 \,dx\right)
\end{equation}
is the Tracy-Widom distribution.
\end{thm}

To state the next result we need the following definitions
\begin{equation}\label{alpha}
\alpha:=\alpha(s)\>=\> \int_{s}^{\iy}q(x)\,u(x)\,dx,
\end{equation}
\begin{equation}\label{mu}
\mu:=\mu(s)\>=\>\int_{s}^{\iy}q(x)dx,
\end{equation}
\begin{equation}\label{nu}
\nu:=\nu(s)\>=\>\int_{s}^{\iy}p(x)dx \>=\> \alpha(s)\>-\> q(s),
\end{equation}

\begin{equation*}
\eta:=\eta(s)=\>\frac{1}{20\sqrt{2}}\int_{s}^{\infty}\hspace{-0.1in}\bigl(6qv+3pu+2p_{2}+
2p_{1}v+2pv_{1}-2q_{2}u-2q_{1}u_{1}-2qu_{2}\bigr)(x)\,dx \,-
\end{equation*}
\begin{equation}
\frac{20c^{2}q'(s)+3p(s)}{20\sqrt{2}}
\end{equation}

\begin{equation*}
E_{c,1}(s)\,=-\frac{1}{20} E_{c,2}(s)\,e^{-\mu }\,
-\,\frac{c\,\alpha}{2 \mu ^2}+\frac{c\, p}{2 \mu }+\frac{(2c-1)\,
\nu ^2}{4 \mu ^2} +c\,u \left(c\,q\,e^{-\mu }\, - \, \frac{\nu }{2
\mu }(1-e^{-\mu } )\right) \>\>+
\end{equation*}
\begin{equation*}
e^{-2 \mu } \left(\frac{\nu\,(\nu+8c\,q) }{32 \mu } -\frac{\eta }{4
\sqrt{2}}\right)\,+ \,e^{-\mu } \left(\frac{2\sqrt{2}\,c^2
q^2-3\,\eta}{4\sqrt{2}}
\,+\,\frac{\nu^{2}-8(2c\,p+c^{2}q^{2})-4c^{2}\alpha^{2}}{32\mu}\right.
\end{equation*}
\begin{equation}\label{Ec1}
-\frac{c^{2}q^{2}}{8\mu^{2}}\,+\,\frac{2-\mu }{2 \mu ^2}
\left(c\,q\,\alpha +\frac{1}{4}\nu^{2}+(c^{2}-c)q^{2}\right)\Biggr)
-\left(4c^{2}\alpha^{2}+3c^{2}q^{2}-\nu^{2}\right)\frac{\cosh(\mu)}{8\mu^{2}}.
\end{equation}

\begin{thm}\cite{Choup3}\label{GOE}
 As $n\ra\iy$
\begin{equation*}
F_{n,1}(\tau(s))^{2}\>=\>F_{2}(s)\cdot\left\{e^{-\mu(s)}\>+\>
\left[c\bigl(q(s)+u(s)\bigr) e^{-\mu(s) } -\frac{\nu(s) }{2 \mu(s)
}\bigl(1-e^{-\mu(s) }\bigr) \right]\,n^{-\frac{1}{3}} \>\> + \right.
\end{equation*}
\begin{equation}\label{Fn1}
 E_{c,1}(s) n^{-\frac{2}{3}}\>\biggr\}\>+\>O(n^{-1})
\end{equation}
uniformly for bounded $s$.
\end{thm}

\subsection{Statement of our results}
If we set
\begin{equation*}
a(t):=\int_{t}^{\iy} q_{n}(x)\,dx \quad \textrm{and} \quad b(t)=\int_{t}^{\iy} p_{n}(x)\,dx
\end{equation*}
then for the Gaussian Orthogonal Ensemble,
\begin{thm}\label{intro close GOE}
For $n$ even,
\begin{equation*}
F_{n,1}(t)^{2}\>=\> F_{n,2}(t)\left\{\frac{1}{2} -c_{\varphi}^{2}\frac{b(t)}{a(t)}  - 2\sqrt{\frac{b(t)}{2a(t)}}c_{\varphi}\sinh\sqrt{2a(t)b(t)}\right.
\end{equation*}
\begin{equation*}
 \left.\left( \frac{1}{2}+c_{\varphi}^{2}\frac{b(t)}{a(t)}\right)\cosh \sqrt{2a(t)b(t)}  \right\}
\end{equation*}
\end{thm}
and for the Gaussian Simplectic Ensemble,

\begin{thm}\label{intro close GSE}
For $n$ odd,
\begin{equation*}
F_{n,4}(t/\sqrt{2})^{2}\>=\> F_{n,2}(t)\frac{1}{2}\left(1+\cosh \sqrt{2a(t)b(t)}\right).
\end{equation*}

\end{thm}
Theorem\eqref{intro close GSE} leads to the following Edgeworth Expansion of $F_{n,4}$

\begin{thm}\label{GSE}

Then as $n\ra\iy$
\begin{equation*}
F_{n,4}(\frac{\tau(s)}{\sqrt{2}})^{2}\>=\>F_{2}(s)\cdot\left\{\cosh^{2}(\frac{\mu}{2})+\frac{c}{2}[
u_{0}(1+\cosh(\mu)) - q \sinh(\mu)] n^{-\frac{1}{3}} \right.
\end{equation*}
\begin{equation}\label{Fn4}
  \frac{1}{4}\left[ \frac{\nu}{2\sqrt{2} \mu} \sinh(\mu)
 +c^{2}q^{2}\cosh(\mu) +\frac{\cosh(\mu) -1}{10} E_{c,2} \right.
 +
\end{equation}
\begin{equation*}
 \left. \sqrt{2}(\eta -\sqrt{2} c^{2}q
u_{0})\sinh(\mu)\biggr]n^{-2/3}
 \right\}
+\frac{\sinh(\mu)}{\mu} O(n^{-1})
\end{equation*}
 uniformly for bounded $s$.
\end{thm}
\vspace{0.3in}

 In \S 2 we complete the derivation of \eqref{f_{n,4}}
as outlined in  \cite{Trac2}. In \S 3 Theorems \ref{intro close GOE} and \ref{intro close GSE} are derived. In \S 4 we justify Theorem \ref{GSE} followed in \S 5 by a brief discussion on the rate of convergence of these distributions as a function of the fine-tuning constant $c$ for large $n$.

\section{Derivation of $F_{n,4}$}
We treat here the case $n$ odd. (Most of the algebra used here can
be found in \cite{Gohb1}, \cite{Gohb2}, \cite{Gohb3}, \cite{Lax1},
and \cite{Whit1}.) We note here that
\begin{equation}\label{repFn4}
F_{n,4}(t/\sqrt{2})^{2}=\det(I-K_{n,4}),
\end{equation}
with

\begin{equation}\label{repKn4}
2K_{n,4}\,=\>\chi \left(\begin{array}{cc}
                K_{n,2}+\psi\otimes \varepsilon \varphi & K_{n,2}\,D\,-\,\psi\otimes \varphi \\
                \varepsilon\,K_{n,2} +\varepsilon\,\psi\otimes \varepsilon\varphi & K_{n,2}+\varepsilon\varphi\otimes  \psi \\
              \end{array}
              \right)\chi.
\end{equation}
We set $J=(t,\iy)$, and $\chi$ represents the multiplication by
the function $\chi_{J}(x)$. This notation allows us to think of the operator with kernel $K_{n,4}$ as an operator acting on $\mathbb{R}$ instead of acting on $(t\>\> \iy)$.\\
We will like to also note that this section is identical to the derivation made in \cite{Trac2}, we only provide here details intentionally left out by the authors of \cite{Trac2}(we suppose to shorten the length of the paper, even if they provided the necessary steps to complete the derivation.)\\

Using the following commutators,
\begin{equation}\label{commutator}
[K_{n,2},D]\>=\> \varphi\otimes\psi +\psi \otimes \varphi,\quad
\quad [\varepsilon, K_{n,2}]\>=\> -\varepsilon \varphi \otimes
\varepsilon \psi -\varepsilon\psi \otimes \varepsilon \varphi
\end{equation}
($\psi$ and $\varphi$ appear as a consequence of the Christophel
Darboux formula applied to $K_{n,2}$,) we have
\begin{equation*}
K_{n,2}+\psi\otimes\varepsilon\varphi=D\,\varepsilon
K_{n,2}+D\,\varepsilon\psi\otimes\varepsilon\varphi=D(\varepsilon
K_{n,2}+\varepsilon\psi\otimes \varepsilon\varphi)=D(
K_{n,2}\,\varepsilon-\varepsilon\varphi\otimes \varepsilon\psi),
\end{equation*}
\begin{equation*}
K_{n,2}\,D-\psi\otimes\varphi=D
\,K_{n,2}+\varphi\otimes\psi=DK_{n,2}+D\varepsilon\varphi\otimes\psi=
D\,(K_{n,2}+\varepsilon\varphi\otimes\psi)
\end{equation*}
and
\begin{equation*}
\varepsilon\,K_{n,2} +\varepsilon\,\psi\otimes
\varepsilon\varphi =K_{n,2}\,\varepsilon
-\varepsilon\,\varphi\otimes \varepsilon\psi
\end{equation*}
as $D\varepsilon\>=\> I$. Our kernel is now
\begin{equation}
2K_{n,4}\,=\> \chi\left(\begin{array}{cc}
                D(
K_{n,2}\,\varepsilon-\varepsilon\varphi\otimes \varepsilon\psi) & D\,(K_{n,2}+\varepsilon\varphi\otimes\psi) \\
                K_{n,2}\,\varepsilon
-\varepsilon\,\varphi\otimes \varepsilon\psi & K_{n,2}+\varepsilon\varphi\otimes\psi \\
              \end{array}
              \right)\chi
\end{equation}
\begin{equation}
=\> \left( \begin{array}{cc} \chi\,D & 0\\ 0& \chi\\
\end{array}\right)\cdot
\left(\begin{array}{cc}
                (
K_{n,2}\,\varepsilon-\varepsilon\varphi\otimes \varepsilon\psi)\chi & (K_{n,2}+\varepsilon\varphi\otimes\psi)\chi \\
                (K_{n,2}\,\varepsilon
-\varepsilon\,\varphi\otimes \varepsilon\psi )\,\chi & (K_{n,2}+\varepsilon\varphi\otimes\psi)\, \chi \\
              \end{array}
              \right)
\end{equation}
Since $K_{n,4}$ is of the form $\mathrm{AB}$, we can use the fact
that $\det(I-\mathrm{AB}) = \det(I-\mathrm{BA})$ and deduce that the
Fredholm determinant of $K_{n,4}$ is unchanged if instead we take
$2K_{n,4}$ to be
\begin{equation}
\left(\begin{array}{cc}
                (
K_{n,2}\,\varepsilon-\varepsilon\varphi\otimes \varepsilon\psi)\chi & (K_{n,2}+\varepsilon\varphi\otimes\psi)\chi \\
                (K_{n,2}\,\varepsilon
-\varepsilon\,\varphi\otimes \varepsilon\psi )\,\chi& (K_{n,2}+\varepsilon\varphi\otimes\psi)\, \chi \\
              \end{array}
              \right)\cdot\left( \begin{array}{cc} \chi\,D & 0\\ 0& \chi\\
\end{array}\right)
\end{equation}
\begin{equation}
\>\>=\>\>\left(\begin{array}{cc}
                (
K_{n,2}\,\varepsilon-\varepsilon\varphi\otimes \varepsilon\psi)\chi\,D & (K_{n,2}+\varepsilon\varphi\otimes\psi)\chi \\
                (K_{n,2}\,\varepsilon
-\varepsilon\,\varphi\otimes \varepsilon\psi )\,\chi\,D & (K_{n,2}+\varepsilon\varphi\otimes\psi)\, \chi \\
              \end{array}
              \right).
\end{equation}
Thus
\begin{equation}\label{fredholm determinant}
\det(I-K_{n,4})\>=\>\det\left(\begin{array}{cc}
               I- \frac{1}{2}(
K_{n,2}\,\varepsilon-\varepsilon\varphi\otimes \varepsilon\psi)\chi\,D & -
\frac{1}{2}(K_{n,2}+\varepsilon\varphi\otimes\psi)\chi \\
                -\frac{1}{2}(K_{n,2}\,\varepsilon
-\varepsilon\,\varphi\otimes \varepsilon\psi )\,\chi\,D & I-\frac{1}{2}(K_{n,2}+\varepsilon\varphi\otimes\psi)\, \chi \\
              \end{array}
              \right).
\end{equation}
Performing row and column operations on the matrix\footnote{This
does not change the determinant, for more details see \cite{Trac2}}
does not change the Fredholm determinant. We first subtract row 1
from row 2, next we add column 2 to column 1 to have the following
matrix
\begin{equation}
\left(\begin{array}{cc}
               I- \frac{1}{2}(
K_{n,2}\,\varepsilon-\varepsilon\varphi\otimes
\varepsilon\psi)\chi\,D -\frac{1}{2}(K_{n,2}+\varepsilon\varphi\otimes\psi)\chi
&\quad -\frac{1}{2}(K_{n,2}+\varepsilon\varphi\otimes\psi)\chi \\
                0 & I \\
              \end{array}
              \right).
\end{equation}
We therefore have,
\begin{equation}\label{fredholm representation}
\det(I-K_{n,4})\>=\> \det\left(I- \frac{1}{2}(
K_{n,2}\,\varepsilon-\varepsilon\varphi\otimes
\varepsilon\psi)\chi\,D -\frac{1}{2}(K_{n,2}+\varepsilon\varphi\otimes\psi)\chi\right)
\end{equation}
\begin{equation}
=\>\det\left(I- K_{n,2}\chi -\frac{1}{2}(
K_{n,2}\,\varepsilon-\varepsilon\varphi\otimes
\varepsilon\psi)\chi\,D +\frac{1}{2}(K_{n,2}-\varepsilon\varphi\otimes\psi)\chi \right)
\end{equation}
\begin{equation*}
=\>\det\left(I- K_{n,2}\chi -\frac{1}{2}(
K_{n,2}\,\varepsilon-\varepsilon\varphi\otimes
\varepsilon\psi)\chi\,D +\frac{1}{2}K_{n,2}\chi+\frac{1}{2}\varepsilon\varphi\otimes\psi\chi -\varepsilon\varphi\otimes\psi\chi \right)
\end{equation*}
\begin{equation*}
=\>\det\left(I- K_{n,2}\chi -\frac{1}{2}
K_{n,2}\,(\varepsilon \chi D-\chi) -\frac{1}{2}\varepsilon\varphi\otimes
\psi (\varepsilon\chi\,D- \chi)  -\varepsilon\varphi\otimes\psi\chi \right)
\end{equation*}
\begin{equation*}
=\>\det\left(I- K_{n,2}\chi -\frac{1}{2}
K_{n,2}\,(\varepsilon \chi D-\varepsilon \,D\,\chi) -\frac{1}{2}\varepsilon\varphi\otimes
\psi (\varepsilon \chi D-\varepsilon \,D\,\chi)  -\varepsilon\varphi\otimes\psi\chi \right)
\end{equation*}
\begin{equation*}
=\>\det\left(I- K_{n,2}\chi -\frac{1}{2}
K_{n,2}\,\varepsilon [\chi\>\> D] -\frac{1}{2}\varepsilon\varphi\otimes
\psi \varepsilon [\chi\>\> D]  -\varepsilon\varphi\otimes\psi\chi \right)
\end{equation*}
We used the fact that $ \varepsilon\, D=I$ to write $\chi=\varepsilon\, D\, \chi$ and the fact that $\varepsilon$ is antisymmetric to have
\begin{equation*}
\varepsilon\varphi\otimes
\varepsilon\psi\,\chi\,D=\varepsilon\varphi\otimes
\psi\varepsilon^{t}\,\chi\,D =-\varepsilon\varphi\otimes
\psi\varepsilon\,\chi\,D.
\end{equation*}
Next we factor out $I-K_{n,2}$ and note that
$(I-K_{n,2})^{-1}=I+R_{n,2}$, where $R_{n,2}$ was defined as the
resolvent of $K_{n,2}$, and
$(I-K_{n,2})^{-1}\varepsilon\varphi=Q_{n,\varepsilon}$. We are
interested on the determinant of the following operator
\begin{equation}\label{product determinants}
\bigl(I-K_{n,2}\chi\bigr)\bigl(I- \frac{1}{2}(K_{n,2}+
R_{n,2}K_{n,2})\varepsilon [\chi\>\> D]-\frac{1}{2}Q_{n,\varepsilon}\otimes
\psi\,\varepsilon [\chi\>\> D]-Q_{n,\varepsilon}\otimes\psi \chi
\bigr).
\end{equation}

 The determinant of the first factor, $\det(I-K_{n,2}\chi)=F_{n,2}$ is a familiar object that was studied in
our work in GUE$_{n}$ see \cite{Choup1} or \cite{Choup2}.
Attention will be given to the second factor of
\eqref{product determinants}. We will represent this factor in the
form $(I-\sum_{j=1}^{k}\alpha_{j}\otimes\beta_{j})$ and use the well
known formula
$\det(I-\sum_{j=1}^{k}\alpha_{j}\otimes\beta_{j})=\det\bigl(\delta_{i,j}-(\alpha_{i},\beta_{j})\bigr)_{i,j=1,\ldots
,k}$ to expand the Fredholm determinant. First we need to find a
representation of $\varepsilon [\chi\,D]$ as a finite rank operator. To
this end we introduce in this section the following notation,
\begin{equation*}
\varepsilon_{k}(x)=\varepsilon(x-a_{k}),\quad
R_{k}(x)=R_{n,2}(x,a_{k}),\quad \delta_{k}(x)=\delta(x-a_{k}),\quad
a_{1}=t,\quad \textrm{and}\quad a_{2}=\iy.
\end{equation*}
With the new notation $J=(t,\iy)=(a_{1},a_{2})$, and the commutator
\begin{equation*}
[\chi\>\>D]=-\delta_{1}\otimes\delta_{1}+\delta_{2}\otimes
\delta_{2},
\end{equation*}
gives
\begin{equation*}
\varepsilon
[\chi\>\>D]=-\varepsilon_{1}\otimes\delta_{1}+\varepsilon_{2}\otimes
\delta_{2}.
\end{equation*}
With this representation, we have
\begin{equation*}
(K_{n,2}+R_{n,2}\,K_{n,2})\varepsilon\,[\chi\>\>D]=\sum_{k=1,2}(-1)^{k}(K_{n,2}+R_{n,2}\,K_{n,2})
\varepsilon_{k}\otimes\delta_{k}
\end{equation*}
and
\begin{equation*}
Q_{n,\varepsilon}\otimes
\psi\,\varepsilon [\chi\>\> D]=\sum_{k=1,2}(-1)^{k}(Q_{n,\varepsilon}\otimes
\psi\,)\cdot(
\varepsilon_{k}\otimes\delta_{k})=\sum_{k=1,2}(-1)^{k}(\psi,\>\>\varepsilon_{k})Q_{n,\varepsilon}\otimes
\delta_{k}.
\end{equation*}

In the last equation we used the formula
$(\alpha\otimes\beta)\cdot(\gamma\otimes \delta)= (\beta,\gamma)
\alpha\otimes \delta$. If we substitute this in the second factor of \eqref{product determinants}, it becomes
\begin{equation}\label{eq4}
I- \frac{1}{2}\sum_{k=1,2}(-1)^{k}(K_{n,2}+
R_{n,2}K_{n,2})\varepsilon_{k}\otimes\delta_{k}-\frac{1}{2}\sum_{k=1,2}(-1)^{k}\bigl(\psi,\varepsilon_{k}\bigr)
Q_{n,\varepsilon}\otimes \delta_{k}
-Q_{n,\varepsilon}\otimes\chi\psi.
\end{equation}
We have
\begin{equation*}
\varepsilon_{2}(x)=\varepsilon_{\infty}(x)=\frac{1}{2}\textrm{sgn}(x-\infty)=-\frac{1}{2},\quad
\textrm{and}\quad R_{2}=R_{n,2}(x,\iy)=0.
\end{equation*}
If we substitute these value in \eqref{eq4}, it then becomes,
\begin{equation}
I-Q_{n,\varepsilon}\otimes \chi \psi
-\frac{1}{2}\bigl[(K_{n,2}+R_{n,2}\,K_{n,2})\varepsilon_{t}+\bigl(\psi,\varepsilon_{t}\bigr)Q_{n,\varepsilon}\bigr]\otimes\delta_{t}+\frac{1}{4}
[(K_{n,2}+R_{n,2}\,K_{n,2})1+\bigl(\psi,1\bigr)Q_{n,\varepsilon}\bigr]\otimes\delta_{\infty}.
\end{equation}
This operator is of the desired form
\begin{equation*}
I-\sum_{k=1,2,3}\alpha_{k}\otimes\beta_{k} \quad \textrm{and} \quad \det(I-\sum_{j=1}^{3}\alpha_{j}\otimes\beta_{j})=\det\bigl(\delta_{i,j}-(\alpha_{i},\beta_{j})\bigr)_{i,j=1,\ldots
,3}
\end{equation*}
with
\begin{equation*}
\alpha_{1}=Q_{n,\varepsilon},\>\>
\end{equation*}
\begin{equation*}
\alpha_{2}=
\frac{1}{2}\bigl[(K_{n,2}+R_{n,2}\,K_{n,2})+ \bigl(\psi,\varepsilon_{t}\bigr)Q_{n,\varepsilon}\bigr],
\end{equation*}
 \begin{equation*}
\alpha_{3}=-\frac{1}{4}
[(K_{n,2}+R_{n,2}\,K_{n,2})1+\bigl(\psi,1\bigr)Q_{n,\varepsilon}\bigr]
\end{equation*}
and
\begin{equation}
\beta_{1}=\chi
\psi,\quad \beta_{2}=\delta_{t} \quad \textrm{and} \quad \beta_{3}=\delta_{\infty}.
\end{equation}
As pointed out in \cite{Trac7}, the contribution of $\beta_{3}$ is zero
\begin{equation*}
(\alpha_{1}\>,\> \beta_{3})=(\alpha_{2}\>,\> \beta_{3})=(\alpha_{3}\>,\> \beta_{3})= 0,
\end{equation*}
thus the determinant reduces to the contribution of $\alpha_{1},\> \alpha_{2},\> \beta_{1}$ and $\beta_{2}$.
The corresponding inner product are;
\begin{equation}
(\alpha_{1},\beta_{1})=\tilde{v}_{n,\varepsilon},\quad
(\alpha_{1},\beta_{2})=q_{n,\varepsilon}+c_{\varphi}=q_{n,\varepsilon} \quad \textrm{as}\quad c_{\varphi}=0 \quad \textrm{for n odd}
\end{equation}
where we have
\begin{equation}
c_{\varphi}=\varepsilon
\,\varphi(\iy)=\frac{1}{2}\int_{-\iy}^{\iy}\varphi(x)\,dx=0 \quad \textrm{as} \quad \varphi \quad \textrm{is odd}
\end{equation}
we also note that
\begin{equation}
c_{\psi}=\varepsilon\,\psi(\iy)=\frac{1}{2}\int_{-\iy}^{\iy}\psi(x)\,dx,
\end{equation}
and for $n$ odd
\begin{equation}\label{c varphi}
c_{\psi}=(\pi\,(n-1))^{1/4}2^{-3/4 -(n-1)/2}\frac{((n-1)!)^{1/2}}{((n-1)/2)!},
\end{equation}
and
\begin{equation}
(\alpha_{2},\beta_{1})=-\frac{1}{2}(\mathcal{P}_{n,4}-a_{4}+a_{4}\tilde{v}_{n,\varepsilon}),\quad
(\alpha_{2},\beta_{2})=-\frac{1}{2}(\mathcal{R}_{n,4}+a_{4}q_{n,\varepsilon}).
\end{equation}
with $a_{4}=(\psi\>,\>\varepsilon_{t})$.
The determinant of \eqref{eq4} is therefore
\begin{equation}\label{determinant}
\bigl(1-\tilde{v}_{n,\varepsilon}(t)\bigr)\bigl(1+\frac{1}{2}\mathcal{R}_{n,4}(t)\bigr)+
       \frac{1}{2}q_{n,\varepsilon}(t)\,\mathcal{P}_{n,4}(t).
\end{equation}

We will make use of large $n$ expansion of the functions appearing \eqref{determinant}  (see \cite{Choup4}) combine with the large $n$ expansion of $F_{n,2}$ to give an edgeworth type expansion of $F_{n,4}$.


\section{Finite $n$ Expression of $F_{n,1}$ and $F_{n,4}$}

In \cite{Choup4} we derived the following results; if we recall that
\begin{equation*}
a(t)=\int_{t}^{\iy}q_{n}(x)\,dx \quad \textrm{and} \quad b(t)=\int_{t}^{\iy}p_{n}(x)\, dx
\end{equation*}
then for the GSE$_{n}$,
\begin{equation}\label{tildevne}
\tilde{v}_{n,\varepsilon}(t)=1-\frac{1}{2}[1+\cosh\sqrt{2a(t)b(t)}]
\end{equation}
and
\begin{equation}\label{qne}
q_{n,\varepsilon}(t)=-\sqrt{\frac{a(t)}{2b(t)}} \sinh\sqrt{2a(t)b(t)}
\end{equation}
for the epsilon independent functions, and the epsilon dependent functions are
\begin{equation}\label{calPn4}
\mathcal{P}_{n,4}(t)=-c_{\psi}\frac{1}{2}[1+\cosh\sqrt{2a(t)b(t)}]-\sqrt{\frac{b(t)}{2a(t)}} \sinh\sqrt{2a(t)b(t)},
\end{equation}
and
\begin{equation}\label{calRn4}
\mathcal{R}_{n,4}(t)=-c_{\psi}\sqrt{\frac{a(t)}{2b(t)}}\sinh\sqrt{2a(t)b(t)}+ \cosh\sqrt{2a(t)b(t)}-1.
\end{equation}
If we recall that
\begin{equation}\label{F_{n,2}}
F_{n,2}(t)\>\>=\>\>\det(I\>-\>K_{n,2})\>\>=\>\> \exp{\biggl(- 2\,
\int_{t}^{\iy}\,(x-t)q_{n}(x)\,p_{n}(x)\,d\,x\biggr)},
\end{equation}

 then have
\begin{thm}\label{close GSE}
For $n$ odd,
\begin{equation*}
F_{n,4}(t/\sqrt{2})^{2}\>=\> F_{n,2}(t)\frac{1}{2}\left(1+\cosh \sqrt{2a(t)b(t)}\right)
\end{equation*}
Or what is the same as 
\begin{equation}\label{close GSE}
F_{n,4}(t/\sqrt{2})\>=\cosh \sqrt{\frac{a(t)b(t)}{2}} \exp{\biggl(- \,
\int_{t}^{\iy}\,(x-t)q_{n}(x)\,p_{n}(x)\,d\,x\biggr)}
\end{equation}
\end{thm}

The surprising fact here is that $c_{\psi}$ drops out and the formula is as simple as possible.  We couldn't wish for a better result.\\

In a similar way we can also use the following representations
\begin{equation*}
F_{n,1}(t)^{2}\>=\>F_{n,2}(t)\left\{(1-\tilde{v}_{n,\varepsilon})(1-\frac{1}{2}\mathcal{R}_{n,1})-\frac{1}{2}(q_{n,\varepsilon}-c_{\varphi})\mathcal{P}_{n,1}\right\}.
\end{equation*}
together with the following representations (see \cite{Choup4} for a derivation of these results),

\begin{equation}\label{tildevne}
\tilde{v}_{n,\varepsilon}(t)=1-\frac{1}{2}[1+\cosh\sqrt{2a(t)b(t)}]+c_{\varphi}\sqrt{\frac{b(t)}{2a(t)}} \sinh\sqrt{2a(t)b(t)}
\end{equation}
and
\begin{equation}\label{qne}
q_{n,\varepsilon}(t)=-\sqrt{\frac{a(t)}{2b(t)}} \sinh\sqrt{2a(t)b(t)} +c_{\varphi} \cosh\sqrt{2a(t)b(t)}.
\end{equation}
We also have

\begin{equation}\label{calPn1}
\mathcal{P}_{n,1}(t)=c_{\varphi}\frac{b(t)}{a(t)}[\cosh\sqrt{2a(t)b(t)}-1]-\sqrt{\frac{b(t)}{2a(t)}} \sinh\sqrt{2a(t)b(t)}
\end{equation}
and

\begin{equation}\label{caltildeRn1}
\mathcal{R}_{n,1}(t)=1-2c_{\varphi}\sqrt{\frac{b(t)}{2a(t)}}\sinh\sqrt{2a(t)b(t)}- \cosh\sqrt{2a(t)b(t)},
\end{equation}

to have  for the GOE$_{n}$,

\begin{thm}\label{close GOE}
For $n$ even,
\begin{equation*}
F_{n,1}^{2}(t)\>=\> F_{n,2}(t)\left\{\frac{1}{2} -c_{\varphi}^{2}\frac{b(t)}{a(t)}  - 2\sqrt{\frac{b(t)}{2a(t)}}c_{\varphi}\sinh\sqrt{2a(t)b(t)}\right.
\end{equation*}
\begin{equation}\label{close GOE}
 \left.\left( \frac{1}{2}+c_{\varphi}^{2}\frac{b(t)}{a(t)}\right)\cosh \sqrt{2a(t)b(t)}  \right\}
\end{equation}
\end{thm}
Unlike the GSE$_{n}$ case, the finite $n$ GOE$_{n}$ does not simplify as well, but it is still a very applicable formula.\\
Theorems \ref{close GOE} and \ref{close GSE} are the main results in this section, they provide a general $n$ formula for the probability distribution function of the largest eigenvalue.  In the next section, we will find a large $n$ expansion of these two functions.

\section{Edgeworth Expansion for GSE$_{n}$}

The previous section tells us that these three probabilities distribution functions $F_{n,1},\>\> F_{n,2}$ and $F_{n,4}$ are all functions of $q_{n}$ and $p_{n}$, thus to find a large $n$ expansion of these functions, we need large $n$ expansion for these functions $q_{n}$ and $p_{n}$. Upon substituting these expansions into the sought after functions \eqref{close GOE} and \eqref{close GSE}, we recover an Edgeworth type expansion of these probability distribution functions. This work was carried out in \cite{Choup1,Choup2} for the GUE$_{n}$, and in \cite{Choup3} for the GOE$_{n}$. In this paper we follow the same technique to find the corresponding expansion for GSE$_{n}$.  This presentation is different from the previous results since it gives a closed form for $F_{n,1}$ and $F_{n,4}$ that can be used to compute these expansions.\\

We quickly see that after substitution of the large $n$ expansion of $q_{n}$ in $a(t)$ and the large $n$ expansion of $p_{n}$ in $b(t)$, the result is

\begin{equation*}
F_{n,1}(\tau(s))\>=\>F_{n,2}(\tau(s))\cdot\left\{e^{-\mu(s)}+\left(\frac{\nu(s)}{2\mu(s)}(e^{-\mu(s)}-1)+cq(s)e^{-\mu(s)}\right)n^{-\frac{1}{3}}\right.
\end{equation*}
\begin{equation*}
+\frac{1}{16 \mu(s)^2}\left(e^{-\mu(s)} \right.\\
\left(\nu(s)^2 \left(\left(-1+e^{\mu(s)}\right) \left(-3+8 c+e^{\mu(s)}\right)+2 \mu(s) \left(1-4 c+4 c^2
\mu(s)\right)\right)-\right.\\
\end{equation*}
\begin{equation*}
8 c \nu(s) \left(-1+e^{\mu(s)}+\mu(s) (-1+2 c \mu(s))\right) \int_{s}^{\iy} q[x] u[x] \, dx+\\
\end{equation*}
\begin{equation*}
8 \mu(s) \left(20 c \left(-1+e^{\mu(s)}\right) \left(\int_{s}^{\iy} q[x] v[x] \, dx-\int_{s}^{\iy} q_1[x] \, dx\right)+\mu(s) \left(\left(-3+20
c^2\right) \right.\right.\\
\end{equation*}
\begin{equation*}
\int_{s}^{\iy} p[x] u[x] \, dx-3 \int_{s}^{\iy} q[x] v[x] \, dx-2 \int_{s}^{\iy} v[x] p_1[x] \, dx-2 \int_{s}^{\iy} p_2[x] \, dx-3 \int_{s}^{\iy} q_1[x] \, dx+\\
\end{equation*}
\begin{equation*}
c^2 \left((\int_{s}^{\iy} q[x] u[x] \, dx)^2-20 \left(2 \int_{s}^{\iy} q[x] u[x]^2 \, dx-3 \int_{s}^{\iy} q[x] v[x] \, dx+\int_{s}^{\iy} q_1[x] \, dx\right)\right)+\\
\end{equation*}
\begin{equation*}
\left.\left.\left.2 \left(\int_{s}^{\iy} u[x] q_2[x] \, dx+\int_{s}^{\iy} q_1[x] u_1[x] \, dx+\int_{s}^{\iy} q[x] u_2[x] \, dx-\int_{s}^{\iy} p[x] v_1[x] \, dx\right)\right)\right)\right)
n^{-\frac{2}{3}}
\end{equation*}
\begin{equation*}
\left. + O(\frac{1}{n}) \right\}
\end{equation*}
\begin{equation}\label{edgeworth GOE}
=\>F_{n,2}(\tau(s))\cdot\left\{e^{-\mu(s)}+\left(\frac{\nu(s)}{2\mu(s)}(e^{-\mu(s)}-1)+cq(s)e^{-\mu(s)}\right)n^{-\frac{1}{3}} +f(s)n^{-\frac{2}{3}} + O(\frac{1}{n}) \right\}
\end{equation}
Using the Edgeworth Expansion of $F_{n,2}$ given in Theorem\ref{Fn2} leads to Theorem\ref{GOE}. For the GSE  we have the following,

\begin{equation*}
F_{n,4}(\frac{\tau(s)}{\sqrt{2}})^{2}\>=\>F_{n,2}(\tau(s))\cdot\left\{ \cosh^{2}(\frac{\mu(s)}{2})-\frac{1}{2}cq(s)\sinh\mu(s)\,n^{-\frac{1}{3}} \right.
\end{equation*}
\begin{equation*}
+\frac{1}{16 \mu(s)}\left(\left(4 c^2 \cosh[\mu(s)] \mu(s) (\nu(s)-\int_{s}^{\iy} q[x] u[x] \, dx)^2+\right.\right.\\
\end{equation*}
\begin{equation*}
\left(-\nu(s)^2+4 \mu(s) \left(\left(3-20 c^2\right) \int_{s}^{\iy} p[x] u[x] \, dx+3 \int_{s}^{\iy} q[x] v[x] \, dx+\right.\right.\\
\end{equation*}
\begin{equation*}
2 \int_{s}^{\iy} v[x] p_1[x] \, dx+2 \int_{s}^{\iy} p_2[x] \, dx+3 \int_{s}^{\iy} q_1[x] \, dx+20 c^2 \left(2 \int_{s}^{\iy} q[x] u[x]^2 \, dx-\right.\\
\end{equation*}
\begin{equation*}
\left.3 \int_{s}^{\iy} q[x] v[x] \, dx+\int_{s}^{\iy} q_1[x] \, dx\right)-2 \left(\int_{s}^{\iy} u[x] q_2[x] \, dx+\int_{s}^{\iy} q_1[x] u_1[x] \, dx+\right.\\
\end{equation*}
\begin{equation*}
\left.\left.\left.\left.\left.\int_{s}^{\iy} q[x] u_2[x] \, dx-\int_{s}^{\iy} p[x] v_1[x] \, dx\right)\right)\right) \sinh[\mu(s)]\right) n^{-\frac{2}{3}} +O(n^{-1})\right\}\\
\end{equation*}
\begin{equation}\label{edgeworth GSE}
=\>F_{n,2}(\tau(s))\cdot\left\{ \cosh^{2}(\frac{\mu(s)}{2})-\frac{1}{2}cq(s)\sinh\mu(s)\,n^{-\frac{1}{3}} + g(s)n^{-\frac{2}{3}} +O(n^{-1})\right\}
\end{equation}

And if we substitute the Theorem\ref{Fn2} we recover Theorem\ref{GSE}

\section{Conclusion}
The major finding about this studies of the probability distribution function of the largest eigenvalues of the classical Gaussian Random Matrix Theory Ensemble, is their dependence on just two quantities, $q_{n}$ and $p_{n}$. Indeed, comparing the representation for these three ensembles, we see that a complete statistical study would only require the knowledge of  of these two functions. The simplicity of representation \eqref{F_{n,2}} and \eqref{close GSE} allows us to have a simple form for the probability density function of the largest eigenvalues, and the somewhat not so simple representation \eqref{close GOE} will do the similar work for the GOE$_{n}$.\\
The second finding is the similarity between the unitary and the Symplectic ensembles, the probability distribution functions are very simple to represent, the fine-tuning constant $c$ does the work, as in these two cases, we can eliminate the first correction term to the Tracy-Widom  limit by setting $c$ to zero. But there are some numerical evidences showing that there is a way to fine-tune the result for the GOE$_{n}$ to speed the convergence rate to the Tracy-Widom distribution. This gives us hope to find the right way to capture this phenomenon for the Orthogonal Ensemble. \\
It would also be convenient to find a simplify version of the functions $f(s)$ appearing in \eqref{edgeworth GOE} and the function $g(s)$ appearing in \eqref{edgeworth GSE}.

\clearpage \vspace{3ex} \noindent\textbf{\large Acknowledgements: }
The author would like to thank Professor Craig Tracy for the
discussions that initiated this work and for the invaluable
guidance, Professor Iain Johnstone for his patience and constant reminder of the importance of this work,and the Department of Mathematical Sciences at the
University of Alabama in Huntsville.


\begin{thebibliography}{10}

\bibitem{Deif1}
J.~Baik, P.~A.~Deift and  K.~Johansson.
\newblock { On the distribution of the length of the
longest increasing subsequence in a random permutation}
\newblock{J. Amer. Math. Soc., 12 (1999), 1119–1178.}

\bibitem{Choup1}
L.~N.~Choup.
\newblock{Edgeworth Expansion of the Largest Eigenvalue Distribution Function
 of GUE and LUE }
\newblock{\em IMRN }Volume 2006, ID 61049, Pages 1-33.

\bibitem{Choup2}
L.~N.~Choup.
\newblock { Edgeworth Expansion of the Largest Eigenvalue Distribution Function of GUE
Revisited}
\newblock {\em J. Math. Phys.} 49, 033508 (2008)

\bibitem{Choup3}
L.~N.~Choup.
\newblock {edgeworth expansion of the largest eigenvalue distribution function of Gaussian orthogonal ensemble }
\newblock {\em J. Math. Phys.} 50, 013512 (2009)

\bibitem{Choup4}
L.~N.~Choup.
\newblock {On Painleve Related Functions Arising in random Matrix Theory }
\newblock {\em to appear}



\bibitem{Deif2}
P.~Deift.
\newblock{Orthogonal Polynomials and Random Matrices: A Riemann-Hilbert Approach}.
\newblock {\em American Mathematical Society}. Courant Lecture
Notes 3, 2000.

\bibitem{Deif3}
P.~Deift,
\newblock{Universality for mathematical and physical systems.}
\newblock{\em International Congress of Mathematicians,} Vol.1, 125-152, Eur.Math.Soc., Z$\ddot{u}$rich, 2007.

\bibitem{Dien1}
M.~Dieng and C.~A.~Tracy.
\newblock{Application of random matrix theory to multivariate statistics}.
\newblock{preprint, Arxiv:math.PR/0603543}.


\bibitem{Fell2}
W.~Feller.
\newblock {\em An Introduction to Probability Theory and Its
Applications} ,Vol.II.
\newblock Second edition, John Wiley, 1971.


\bibitem{Gohb1}
I.~Gohberg, S.~Goldberg, and M.~A. Kaashoek.
\newblock {\em {Classes of Linear Operators, Vol. I}}, volume~49 of {\em
  Operator Theory: Advances and Applications}.
\newblock Birkh{\"a}user, 1990.

\bibitem{Gohb2}
I.~Gohberg, S.~Goldberg, and M.~A. Kaashoek.
\newblock {\em {Classes of Linear Operators, Vol. II}}, volume~63 of {\em
  Operator Theory: Advances and Applications}.
\newblock Birkh{\"a}user, 1993.

\bibitem{Gohb3}
I.~C. Gohberg, M.~G. Kre$\breve{i}$n.
\newblock {\em {Introduction to the Theory of Linear Nonselfadjoint Operators}}, volume~18 of {\em
  Translations of Mathematical Monographs}.
\newblock American Mathematical Society, 1969.



\bibitem{Hoch1}
H.~Hochstadt.
\newblock {\em {The Functions of Mathematical Physics}}, volume~23 of {\em
  Pure and Applied Mathematics: A series of texts and Monographs }.
\newblock Wiley-Interscience, 1971.

\bibitem{Johansson}
K.~Johansson.
\newblock{ Toeplitz determinants, random growth and determinantal
processes.}
\newblock{\em Proceedings of the ICM, Beijing 2002,} vol. 3,
53--62, math.PR/0304368.

\bibitem{John1}
I.~M.~Johnstone,
\newblock {On the distribution of the largest eigenvalue in principal component
  analysis},
\newblock {\em Ann. Stats.}, 29(2):295--327, 2001.


\bibitem{Lax1}
P.~D.~Lax.
\newblock{ \em Functional Analysis}
\newblock{ Wiley-Interscience}, 2002.


\bibitem{Meht1}
M.~L.~Mehta.
\newblock {\em Random Matrices, Revised and Enlarged Second Edition}.
\newblock Academic Press, 1991.

%
%

\bibitem{Szeg1}
G.~Szeg\"o.
\newblock{Orthogonal Polynomials.}
\newblock American Mathematical Society Colloquium Publications
Volume 23

\bibitem{Taylor}
M.~E.~Taylor
\newblock{Partial Differential Equations}.
\newblock{Springer-Verlag, New York, 1996}

\bibitem{Trac3}
C.~A.~Tracy and H.~Widom.
\newblock {Level--spacing distributions and the Airy kernel}.
\newblock {\em Commun. Math. Physics}, 159:151--174, 1994.

\bibitem{Trac7}
C.~A.~Tracy and H.~Widom.
\newblock {Fredholm determinants, differential equations and matrix models}.
\newblock {\em Commun. Math. Physics}, 163:33--72, 1994.

\bibitem{Trac2}
C.~A.~Tracy and H.~Widom.
\newblock {On orthogonal and symplectic matrix ensembles}.
\newblock {\em Commun. Math. Physics}, 177:727--754, 1996.

\bibitem{Trac1}
C.~A.~Tracy and H.~Widom.
\newblock {Correlation functions, cluster functions, and spacing distributions
  for random matrices}.
\newblock {\em J. Stat. Phys.}, 92(5--6):809--835, 1998.

\bibitem{Trac4}
C.~A. Tracy and H.~Widom.
\newblock {Airy kernel and Painlev\'e II}.
\newblock In {\em Isomonodromic deformations and applications in physics},
  volume~31 of {\em {CRM Proceedings \& Lecture Notes}}, pages 85--98. Amer.
  Math. Soc., Providence, RI, 2002.

\bibitem{Trac8}
C.~A.~Tracy and H.~Widom.
\newblock {Distribution functions for largest eigenvalues
and their applications}.
\newblock In {\em Proceedings of the International Congress
of Mathematicians, Beijing 2002}, Vol.~I, ed. LI Tatsien, Higher
Education Press, Beijing, pgs.~587--596, 2002.

\bibitem{Trac5}
C.~A.~Tracy and H.~Widom.
\newblock {Matrix kernels for the Gaussian orthogonal and symplectic
  ensembles}.
\newblock {\em Ann. Inst. Fourier, Grenoble}, 55, 2197--2207, 2005.

\bibitem{Whit1}
E.~T.~Whittaker and G.~N.~Watson.
\newblock{\em A Course of Modern Analysis} Fourth Edition
\newblock{Cambridge University Press}, 2004.


\end{thebibliography}
\end{document}